\documentclass[12pt]{article}

\usepackage[margin=1in]{geometry}
\usepackage{amsmath}
\usepackage{bm}
\usepackage{url}
\usepackage{xcolor}
\usepackage{graphicx}
\usepackage{booktabs}
\usepackage{siunitx}
\usepackage{amssymb} 
\usepackage[numbers,sort&compress]{natbib} 
\usepackage{hyperref} 

\usepackage{subcaption}

\newcommand{\e}[1]{\ensuremath{\times 10^{#1}}}

\title{A $C^1$-continuous finite element formulation for solving the Jeffery-Hamel boundary value problem}
\author{J.~W.~Peterson \and R.~H.~Stogner}
\date{\today}

\begin{document}
\maketitle

\begin{abstract}
The third-order Jeffery-Hamel ODE governing the flow of an incompressible
fluid in a two-dimensional wedge is briefly derived, and a $C^1$ finite element formulation of
the equation is developed.  This formulation has several advantages,
including a natural framework for enforcing the boundary conditions, a
numerically efficient solution procedure, and suitability for
implementation within well-established, open, scientific computing
tools.  The finite element formulation is shown to be non-coercive,
and therefore not ideal for proving existence, uniqueness, or \emph{a
priori} error estimates, but the numerical solutions computed with
quartic Hermite elements are nevertheless found to converge to
reference solutions at nearly optimal rates ($\mathcal{O}(h^4)$ in
both $L^2$ and $H^1$ norms).  Further work is required to better
understand the cause of the suboptimal convergence rates, and a linear
model problem which exhibits analogous characteristics is also
discussed as a possible starting point for future theoretical
analyses.
\end{abstract}

\section{Introduction}
Viscous, incompressible flow in a two-dimensional wedge, frequently
referred to as Jeffery-Hamel flow, is described in many references
dating back to the original works by Jeffery~\cite{Jeffery_1915} and
Hamel~\cite{Hamel_1916}, the comprehensive discussion of the
various possible configurations of the flow field and a general
solution method in terms of elliptic functions by
Rosenhead~\cite{Rosenhead_1940}, and the modern treatment in fluid
mechanics textbooks~\cite{White_1991,Panton_2006,Langlois_2014}.
Recent research~\cite{Haines_2011} has focused on performing
nonlinear stability analyses, computing bifurcation diagrams, and classifying
non-unique, stable solutions in various parameter regimes.

In addition to fundamental fluid mechanics research, the Jeffery-Hamel
flow solutions are also of great utility as validation tools for
finite difference, finite element, and related numerical codes
designed to solve the incompressible Navier-Stokes equations under
more general conditions. Although a closed-form analytical solution to
the Jeffery-Hamel equations is known, it is typically more convenient
to work with the solution in a ``semi-analytical'' form, that is, a
form which can be obtained to any desired accuracy using (yet another)
numerical method.  In this short note, we describe a numerical method
based on a $C^1$ finite element formulation for efficiently and
accurately approximating solutions to the Jeffery-Hamel equations, and
compare it to other established techniques in terms of computational
expense and implementation difficulty.

The rest of the paper is arranged as follows: the governing equations
and the assumptions leading up to them are described
in~\S\ref{sec:goveqn}.  The classical procedure which leads to the
semi-analytical form of the solution is described
in~\S\ref{sec:semianalytical}.  In~\S\ref{sec:fem}, the $C^1$ finite
element method is described. Theoretical aspects related to the
existence and uniqueness of solutions to the finite element
formulation are discussed in~\S\ref{sec:exuniq}.  The numerical
results are presented in~\S\ref{sec:results}, and their accuracy is
compared to established methods.  Finally, in~\S\ref{sec:future},
we summarize the conclusions of this research, and some
directions for future work are listed and described.

\section{Governing equations\label{sec:goveqn}}
The Jeffery-Hamel solution corresponds to flow
constrained to a wedge-shaped region:
$r_1 \leq r \leq r_2$, $-\alpha \leq \theta \leq \alpha$.
The origin ($r=0$) is a singular point of the flow, and is always
excluded from numerical computations.  The governing equations are the
incompressible Navier-Stokes mass and momentum conservation equations in
cylindrical polar coordinates.  It is assumed that the flow is purely
radial ($u_{\theta} = 0$), and the boundary conditions are no slip on
the solid walls ($u_r(r, \pm\alpha)=0$) and symmetry
about the centerline ($\theta=0$) of the channel.  Under these
assumptions, the incompressible Navier-Stokes equations simplify to:
\begin{align}
  \label{eqn:mass}
  \frac{1}{r}\frac{\partial (r u_r)}{\partial r} &= 0
  \\
  \label{eqn:rmomentum}
  u_r \frac{\partial u_r}{\partial r} &= -\frac{1}{\rho}\frac{\partial p}{\partial r}
                                        + \nu\left( \frac{\partial^2 u_r}{\partial r^2} +
                                        \frac{1}{r}\frac{\partial u_r}{\partial r} -
                                        \frac{u_r}{r^2} +
                                        \frac{1}{r^2} \frac{\partial^2 u_r}{\partial \theta^2}\right)
  \\
  \label{eqn:thetamomentum}
  0 &= -\frac{1}{\rho}\frac{\partial p}{\partial \theta} + \frac{2\nu}{r}\frac{\partial u_r}{\partial \theta}
\end{align}
where the dynamic viscosity, $\mu$, density $\rho$, and kinematic
viscosity $\nu \equiv \frac{\mu}{\rho}$ are given constants which
depend on the fluid.

\section{Semi-analytical solution strategy\label{sec:semianalytical}}
The solution to~\eqref{eqn:mass} is particularly simple, and inspires
the non-dimensionalization and the form of the eventual solution to
the problem.  Integrating~\eqref{eqn:mass} gives:
\begin{align}
  \label{eqn:mass_integrated}
  r u_r = F(\theta)
\end{align}
where $F(\theta)$ is a function that depends
only on the angular coordinate, $\theta$.
Eqn.~\eqref{eqn:mass_integrated} states that the quantity $r u_r$ is
\emph{constant} along any fixed angular direction $\theta=\text{const}$.
Since we expect the maximum velocity to occur along the centerline (due to the no-slip boundary conditions on the solid walls),
we define the quantity $u_{\text{max}} (r) \equiv u_r(r,0)$.  This lets us define the constant
\begin{align}
  \lambda \equiv r u_{\text{max}} (r) = \text{const}
\end{align}
That is, the centerline velocity varies with $r$ throughout the
domain, but the product $r u_{\text{max}} (r)$ remains fixed.  The
quantity $\lambda$ has physical units of $\frac{L^2}{T}$, and
allows us to define the dimensionless Reynolds number for this problem as
\begin{align}
  \text{Re} \equiv \frac{\lambda \alpha}{\nu}
\end{align}
The quantity $\lambda$ can always be computed once
Re, $\alpha$, and the fluid property $\nu$ have been specified.  Finally, normalizing the angular
coordinate according to $\eta \equiv \frac{\theta}{\alpha}$, we obtain the
non-dimensional form of~\eqref{eqn:mass_integrated} as
\begin{align}
  \label{eqn:nondimensional_f}
  \frac{r u_r}{\lambda} = f(\eta)
\end{align}
where $f$ is still an unknown, dimensionless function that depends
only on $\eta$, and must satisfy several boundary conditions to be described later.
We can rearrange~\eqref{eqn:nondimensional_f} as:
\begin{align}
  \label{eqn:ur}
  u_r = \frac{\lambda }{r} f(\eta)
\end{align}
and then, after making the following substitutions
\begin{align}
  \frac{\partial u_r}{\partial r}
  &=
  -\frac{\lambda}{r^2}f
  \\
  \frac{\partial^2 u_r}{\partial r^2}
  &= \frac{2\lambda}{r^3}f
  \\
  \frac{\partial u_r}{\partial \theta}
  &= \frac{\lambda}{r} \frac{\partial f}{\partial \eta} \frac{\partial \eta}{\partial \theta}
  \equiv \frac{\lambda}{r\alpha} f'
  \\
  \label{eqn:ur_theta_theta}
  \frac{\partial^2 u_r}{\partial \theta^2}
  &= \frac{\lambda}{r\alpha^2} f''
\end{align}
in~\eqref{eqn:rmomentum} and~\eqref{eqn:thetamomentum}, we obtain:
\begin{align}
  \label{eqn:rmomentum_nondim}
  -\frac{\lambda^2}{r^3}f^2 &= -\frac{1}{\rho}\frac{\partial p}{\partial r}
                                        + \frac{\nu \lambda}{r^3 \alpha^2} f''
  \\
  \label{eqn:thetamomentum_nondim}
  0 &= -\frac{1}{\rho \alpha}\frac{\partial p}{\partial \eta} + \frac{2\nu\lambda}{r^2 \alpha}f'
\end{align}
Multiplying~\eqref{eqn:rmomentum_nondim} and~\eqref{eqn:thetamomentum_nondim} by $\rho$ and rearranging gives:
\begin{align}
  \label{eqn:rmomentum_nondim_rearranged}
  \frac{\partial p}{\partial r} &= \frac{\lambda \mu}{\alpha^2 r^3}\left(f'' + \alpha\text{Re}f^2\right)
  \\
  \label{eqn:thetamomentum_nondim_rearranged}
  \frac{\partial p}{\partial \eta} &= \frac{2\mu\lambda}{r^2}f'
\end{align}
The next step is to eliminate $p$
from~\eqref{eqn:rmomentum_nondim_rearranged}
and~\eqref{eqn:thetamomentum_nondim_rearranged} in order to solve for
$f$.  This is accomplished by
differentiating~\eqref{eqn:rmomentum_nondim_rearranged} with respect
to $\eta$, and~\eqref{eqn:thetamomentum_nondim_rearranged} with
respect to $r$, and subtracting.  The result is:
\begin{align}
  \label{eqn:ode3}
  f''' + 2\text{Re}\,\alpha f\!f' + 4\alpha^2f' = 0
\end{align}
Equation~\eqref{eqn:ode3} is a third-order boundary value problem whose description is completed by
the specification of the following three boundary conditions:
\begin{align}
  \label{eqn:f0}
  f(0)  &= 1 \quad \text{(centerline velocity)} \\
  f'(0) &= 0 \quad \text{(centerline symmetry)} \\
  \label{eqn:f1}
  f(1)  &= 0 \quad \text{(no slip)}
\end{align}
Equation~\eqref{eqn:ode3} has an analytical solution which is given in terms of
elliptic integrals, but it is more common (and in many respects simpler) to instead compute
a highly-accurate approximate solution to~\eqref{eqn:ode3} using a numerical method.
One possible numerical approach is to rewrite~\eqref{eqn:ode3} as a system of three
first-order ODEs, and use a ``shooting method'' to iteratively compute
solutions until the initial data which produces the desired end
condition at $\eta=1$ is obtained.  Another possibility is to
solve~\eqref{eqn:ode3} as a boundary value problem using any of a
number of numerical procedures which have been developed for this
class of problem.  The $C^1$ finite element solution pursued in the
present work falls into this category, and is discussed in further detail in \S\ref{sec:fem}.

Once $f$ has been computed numerically, $u_r$ follows directly but there is an additional step required to find $p$.
Integrating~\eqref{eqn:rmomentum_nondim_rearranged} with respect to $r$ yields
\begin{align}
  \label{eqn:rmomentum_integrated}
  p = p^{\ast} - \frac{\lambda \mu}{2 \alpha^2 r^2}\left(f'' + \alpha\text{Re}f^2\right) + T(\theta)
\end{align}
where $p^{\ast}$ is an arbitrary constant, and $T(\theta)$ is a function of $\theta$ only.
Similarly, integrating~\eqref{eqn:thetamomentum_nondim_rearranged} with respect to $\eta$ gives:
\begin{align}
  \label{eqn:thetamomentum_integrated}
  p = p^{\ast} + \frac{2\mu \lambda}{r^2} f + R(r)
\end{align}
where $R(r)$ is a function of $r$ only.  If we make the particular choices
\begin{align}
  T(\theta) &\equiv 0
  \\
  \label{eqn:R}
  R(r) &\equiv \frac{2\mu\lambda K}{r^2}
\end{align}
where $K$ is a constant, then the pressure fields defined by~\eqref{eqn:rmomentum_integrated} and~\eqref{eqn:thetamomentum_integrated}
are the same if and only if:
\begin{align}
  - \frac{\lambda \mu}{2 \alpha^2 r^2}\left(f'' + \alpha\text{Re}f^2\right) = \frac{2\mu \lambda}{r^2} (f+K)
\end{align}
or, solving for $K$ in terms of the (now) known function $f$:
\begin{align}
  \label{eqn:K}
  K = -\frac{1}{4 \alpha^2}\left(f'' + \alpha\text{Re}f^2\right) -f
\end{align}
Multiplying~\eqref{eqn:K} by $f'$ and integrating from 0 to 1 gives:
\begin{align}
  \label{eqn:K_integrated}
  K \int_0^1 f' \text{d}\eta = -\frac{1}{4\alpha^2} \int_0^1 \left[f' \left(f'' + \alpha\text{Re} f^2\right) + 4\alpha^2 ff' \right]\text{d}\eta
\end{align}
Applying integration by parts to~\eqref{eqn:K_integrated} then results in
\begin{align}
  \label{eqn:K_integrated_2}
  K \left. f \right|_0^1 = -\frac{1}{4\alpha^2} \left[ \frac{1}{2} \left. (f')^2 \right|_0^1 + \frac{\alpha \text{Re}}{3} \left. f^3 \right|_0^1 + 2\alpha^2 \left. f^2 \right|_0^1 \right]
\end{align}
Finally, substituting in the boundary conditions~\eqref{eqn:f0}--\eqref{eqn:f1} gives
\begin{align}
  \label{eqn:K_integrated_3}
  K  = \frac{1}{4\alpha^2} \left(\frac{1}{2} f'(1)^2 - \frac{\alpha \text{Re}}{3} - 2\alpha^2\right)
\end{align}
that is, the value of $K$ depends only on the known constants of the problem and the gradient $f'(1)$ at the right-hand boundary.
Once $f$, and consequently $K$, are known, the pressure is given (up to an arbitrary constant) by:
\begin{align}
  \label{eqn:pressure_final}
  p = p^{\ast} + \frac{2\mu \lambda}{r^2} (f + K)
\end{align}
In a numerical simulation, the arbitrary pressure constant can be
selected by ``pinning'' a single value of the pressure wherever it is
convenient, typically on the boundary.  Numerically computed $K$
values for some representative (Re, $\alpha$) values are given in
Table~\ref{tab:K}.
\begin{table}[hbt]
  \centering
  \caption{Values of $K$ for some representative $(\text{Re},\alpha)$
    combinations.  When the product $\text{Re}\,\alpha > 0$, there is an
    increasingly strong adverse pressure gradient (represented by larger
    negative $K$ values), and for $\text{Re}\,\alpha > 10.31$ the profiles
    are linearly unstable~\cite{Rosenhead_1940}.  When the product
    $\text{Re}\,\alpha < 0$ (hence $\lambda < 0$), there is ``converging''
    flow (toward the origin) and the solutions are linearly
    stable.\label{tab:K}}
  \vspace{4pt}
  \begin{tabular}{ll}
    \toprule
    $(\text{Re},\alpha)$     & $K$ \\
    \midrule
    $(30,15^{\circ})$   & $-9.7822146449$                \\ 
    $(110,3^{\circ})$   & $-1.4387160807\e{2}$           \\ 
    $(-80,5^{\circ})$   & $\phantom{-}2.5439853775\e{2}$ \\ 
    \bottomrule
  \end{tabular}
\end{table}

\section{Finite element formulation\label{sec:fem}}
The numerical solution component of the Jeffery-Hamel equations has
been tackled by a wide variety of approximation methods over the
years.  The reasons for the popularity of the equations are quite varied, and include
their utility as code verification tools, the ease with which results
can be verified against tabulated values in the literature, and the
interesting mathematical characteristics---including nonlinearity and
higher derivatives---possessed by the equations.

Solution techniques include boundary value problem
solvers~\cite{Shampine_2006,Boisvert_2012}, the modified decomposition
method~\cite{Kezzar_2015a,Kezzar_2015b}, the reproducing kernel
Hilbert space method~\cite{Inc_2013}, homotopy
methods~\cite{Motsa_2010,Esmaeilpour_2010,Joneidi_2010}, integral
transform methods~\cite{Sushila_2013}, and mixed analytical/numerical
solution methods based on computer algebra
software~\cite{Corless_2007}.  In this work, we pursue a $C^1$ finite
element solution of~\eqref{eqn:ode3} in order to show that this
variational approach is capable of achieving accurate results in a
computationally efficient manner.  The prevalence of open source,
customizable, and high-quality finite element
libraries~\cite{Kirk_2006,Bangerth_2007,Logg_2012} greatly simplifies
the task of implementing such solution algorithms, and helps ensure
correct code and the propagation of reproducible, curated results.

The finite element method proceeds by multiplying~\eqref{eqn:ode3} by
a test function $v \in H^2(\Omega)$, the Hilbert space of functions
with square-integrable second derivatives on $\Omega = (0,1)$, and
integrating over the domain to obtain:
\begin{align}
  \label{eqn:ode3_weak}
  \int_0^1 \left( f''' + 2\text{Re}\,\alpha f\!f' + 4\alpha^2f' \right) v \; \text{d}x = 0
\end{align}
Integrating by parts twice on the first term produces:
\begin{align}
  \label{eqn:ode3_weak_ibp}
  \int_0^1 f'\left( v'' + 2\text{Re}\,\alpha f v + 4\alpha^2 v \right) \text{d}x + \left. f''v \right|_0^1 - \left. f'v' \right|_0^1 = 0
\end{align}
We then incorporate the boundary conditions~\eqref{eqn:f0}--\eqref{eqn:f1} into the test and trial spaces by defining
\begin{align}
  V_0 &= \{v : v \in H^2(\Omega), v(0) = v'(0) = v(1) = 0\} \subset H^2(\Omega)
  \\
  S &= \{u : u \in H^2(\Omega), u(0) = 1, u'(0) = u(1) = 0\} \subset H^2(\Omega)
\end{align}
and seek $f \in S$ satisfying:
\begin{align}
  \label{eqn:ode3_weak_ibp_with_bcs}
  \int_0^1 f'\left( v'' + 2\text{Re}\,\alpha f v + 4\alpha^2 v \right) \text{d}x - f'(1)v'(1) = 0 \quad \forall\, v \in V_0
\end{align}
Introducing a mesh and the finite-dimensional subspaces $V_0^h \subset V_0$, $S^h \subset S$
spanned by the basis $\{\phi_i\}, i=1..N$ leads to the discrete residual statement: find $f_h \in S^h$ such that
$ R_i(f_h) = 0$ for $i=1..N$, where
\begin{align}
  \label{eqn:Ri}
  R_i(f_h) \equiv \int_0^1 f_h'\left( \phi_i'' + 2\text{Re}\,\alpha f_h\phi_i + 4\alpha^2 \phi_i \right) \text{d}x - f_h'(1)\phi_i'(1)
\end{align}
The associated Jacobian contribution is given by:
\begin{align}
  \label{eqn:Jij}
  J_{ij}(f_h) \equiv \int_0^1 \left[
  f_h'\left( 2\text{Re}\,\alpha \phi_j\phi_i \right) +
  \phi_j'\left( \phi_i'' + 2\text{Re}\,\alpha f_h\phi_i + 4\alpha^2 \phi_i \right)
  \right] \text{d}x - \phi_j'(1)\phi_i'(1)
\end{align}
The nonlinear system of equations defined by~\eqref{eqn:Ri} can
be solved for $f_h$ using e.g.\ an inexact Newton method which employs
high-performance sparse preconditioned Krylov solvers at each
iteration.

\begin{figure}[hbt]
  \centering
  \includegraphics[width=.7\linewidth]{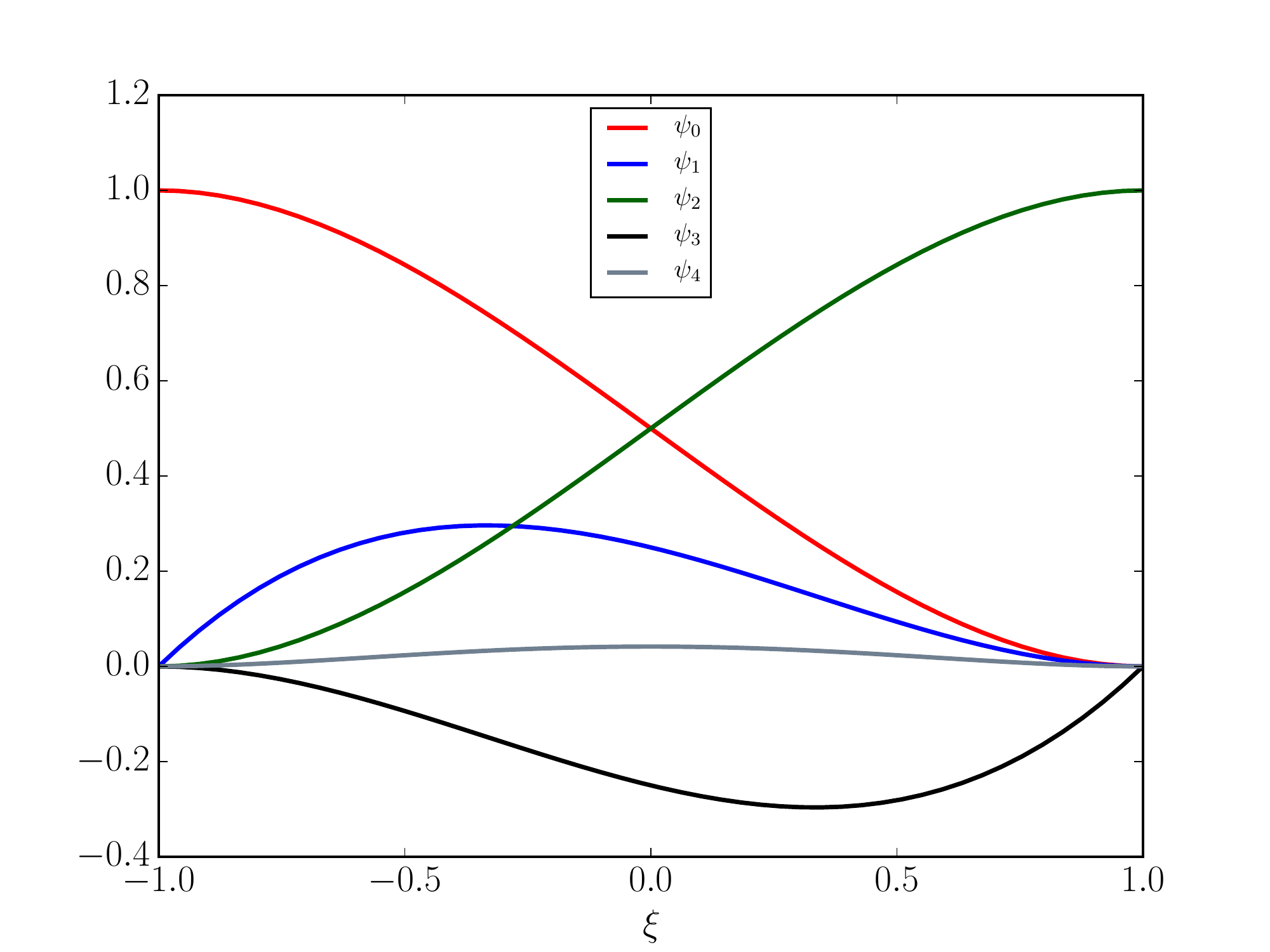}
  \caption{Hermite basis functions up to order $p=4$ in 1D.\label{fig:hermite}}
\end{figure}

The simplest and most natural family of finite element shape functions
which satisfies the requirements of~\eqref{eqn:Ri} are the Hermite
elements, which are composed of $C^1$-continuous polynomials for any
order $p\geq 3$.  The first four element shape functions
(shown in Fig.~\ref{fig:hermite}) correspond to the value and gradient
degrees of freedom at the left and right nodes, while the higher-order
basis functions are ``bubbles.'' The degrees of freedom associated to
the bubble functions could be statically-condensed out of the linear
systems before solution, but we do not pursue this optimization in the
present work.  Finally, we note that the residual~\eqref{eqn:Ri} and
Jacobian~\eqref{eqn:Jij} contributions require a quadrature rule
capable of evaluating polynomials of order $3p-1$ exactly when the
underlying basis is of order $p$.  For $p=3$, this corresponds to a
five point Gauss quadrature rule, while for $p=4$, a six point rule is
required.

\section{Existence and uniqueness of solutions\label{sec:exuniq}}
It is reasonable to ask whether a solution to the nonlinear
variational problem~\eqref{eqn:ode3_weak_ibp_with_bcs} exists, and if
so, whether it is unique. One related theorem is discussed
in~\cite{Zeidler_1990}, pg.~472. In abstract notation, the nonlinear
problem
\begin{align}
  A(f) = 0, \quad f \in X
\end{align}
where $X$ is a reflexive Banach space, has a solution if
\begin{enumerate}
\item The operator $A : X \rightarrow X^{\ast} $ is monotone, i.e.
  \begin{align}
    \label{eqn:monotonicity_def}
    \langle A(f) - A(g), f - g \rangle \geq 0 \quad \forall\, f, g \in X
  \end{align}
\item $A$ is hemicontinuous, i.e. $A(f + tg)$ converges weakly to $A(f)$ $\forall f,g \in X$ as $t \rightarrow 0^+$.
\item $A$ is coercive, i.e.
  %
  \begin{align}
    \lim_{\|f\|_{X} \rightarrow \infty} \frac{\langle A(f), f \rangle}{\|f\|_{X}} = \infty
  \end{align}
\end{enumerate}
where the duality pairing is defined in terms of~\eqref{eqn:ode3_weak_ibp} as
\begin{align}
  \label{eqn:duality_pairing}
  \langle A(f), v\rangle \equiv
  \int_0^1 f'\left( v'' + 2\text{Re}\,\alpha f v + 4\alpha^2 v \right) \text{d}x + \left. f''v \right|_0^1 - \left. f'v' \right|_0^1
\end{align}
Unfortunately, it is easy to see that $A(f)$ does not satisfy preconditions 1 and 3 above.
For example, to show that $A(f)$ is not coercive, we can directly compute
\begin{align}
  \nonumber
  \langle A(f), f\rangle &=
  \int_0^1 f'\left( f'' + 2\text{Re}\,\alpha f^2 + 4\alpha^2 f \right) \text{d}x + \left. f''f \right|_0^1 - \left. (f')^2 \right|_0^1
  \\
  \nonumber
  &= \int_0^1 \left( \frac{1}{2} \left[(f')^2\right]' + \frac{2\text{Re}\,\alpha}{3} (f^3)' + 2\alpha^2 (f^2)' \right) \text{d}x
    + \left. f''f \right|_0^1 - \left. (f')^2 \right|_0^1
  \\
  \nonumber
  &= -\frac{1}{2} \left. (f')^2 \right|_0^1 + \frac{2\text{Re}\,\alpha}{3} \left. f^3 \right|_0^1 + 2\alpha^2 \left. f^2 \right|_0^1 + \left. f''f \right|_0^1
  \\
  \label{eqn:duality_pairing_ff}
  &= -\frac{1}{2} f'(1)^2 - \frac{2\text{Re}\,\alpha}{3}  - 2\alpha^2  - f''(0)
\end{align}
where the last line follows by imposing the boundary conditions~\eqref{eqn:f0}--\eqref{eqn:f1}.
Thus $\langle A(f), f\rangle$ is not bounded from below by any multiple of $\|f\|_{H^2}$, and
we conclude that the duality pairing~\eqref{eqn:duality_pairing} is not coercive.

To show that $A(f)$ is not monotone, let $q\equiv f-g$ for brevity, and directly compute:
\begin{align}
  \nonumber
  \langle A(f) - A(g), q \rangle &= \int_0^1 \left( q' q'' + 2\text{Re}\,\alpha (f-g) (ff' - gg') + 4\alpha^2 q q' \right) \text{d}x + \left.q''q\right|_0^1 - \left. (q')^2 \right|_0^1
  \\
  \label{eqn:monotonicity_show}
  &= 2\text{Re}\,\alpha \int_0^1  (f-g) (ff' - gg') \; \text{d}x
  + \left.q''q\right|_0^1 - \frac{1}{2} \left. (q')^2 \right|_0^1 + \left. 2\alpha^2q^2\right|_0^1
\end{align}
To show lack of monotonicity, we need only find a single $f^{\ast}$ and $g^{\ast}$ for which~\eqref{eqn:monotonicity_def}
does not hold.  For simplicity, assume that $f^{\ast}=g^{\ast}$ on the boundary, and therefore the boundary terms vanish
in~\eqref{eqn:monotonicity_show} vanish, leaving
\begin{align}
  \label{eqn:monotonicity_star}
  \langle A(f^{\ast}) - A(g^{\ast}), f^{\ast} - g^{\ast} \rangle
  &= 2\text{Re}\,\alpha \int_0^1  (f^{\ast}-g^{\ast}) (f^{\ast} f^{\ast\prime} - g^{\ast}g^{\ast\prime}) \; \text{d}x
\end{align}
Next, assume that for this specific choice of $f^{\ast}$ and $g^{\ast}$, the operator is strictly monotone, i.e.
\begin{align}
  \langle A(f^{\ast}) - A(g^{\ast}), f^{\ast} - g^{\ast} \rangle > 0
\end{align}
Letting $u = -f^{\ast}$ and $v = -g^{\ast}$ in~\eqref{eqn:monotonicity_star} then gives
\begin{align}
  \nonumber
  \langle A(u) - A(v), u - v \rangle
  &= -2\text{Re}\,\alpha \int_0^1  (f^{\ast}-g^{\ast}) (f^{\ast} f^{\ast\prime} - g^{\ast}g^{\ast\prime}) \; \text{d}x
  \\
  \label{eqn:monotonicity_star2}
  & < 0
\end{align}
and therefore the operator $A(f)$ is not monotone.

An existence and uniqueness proof is possible if we instead formulate
the problem as a system of nonlinear first-order ODEs by defining: $y_0 \equiv
f$, $y_1 \equiv f'$, and $y_2 \equiv f''$.  The third-order ODE~\eqref{eqn:ode3} can
then be written as
\begin{align}
  \label{eqn:y0}
  y_0' &= y_1
  \\
  y_1' &= y_2
  \\
  \label{eqn:y2}
  y_2' &= -2\text{Re}\,\alpha y_0 y_1 - 4\alpha^2 y_1
\end{align}
for $\eta \in [0, 1]$ (since $f(\eta)$ is symmetric about $\eta=0$)
subject to the initial conditions:
\begin{align}
  \label{eqn:y0_ic}
  y_0(0) &= 1
  \\
  \label{eqn:y1_ic}
  y_1(0) &= 0
  \\
  \label{eqn:y2_ic}
  y_2(0) &= s
\end{align}
where $s$ is unknown, and must be determined iteratively to ensure
that the end condition $y_0(1) = 0$ is satisfied, for example via the
shooting method. We can then write equations~\eqref{eqn:y0}--\eqref{eqn:y2},
\eqref{eqn:y0_ic}--\eqref{eqn:y2_ic} as
\begin{align}
  \label{eqn:ivp}
  \vec{y}^{\,\prime} &= \vec{F}(\vec{y})
  \\
  \vec{y}(0) &= \vec{y}_0
\end{align}
By the Picard-Lindel\"{o}f theorem~\cite{Coddington_1955}, if
$\vec{F}(\vec{y})$ is Lipschitz continuous in $B_r(\vec{y}_0)$ (a
closed ball of radius $r$ centered at $\vec{y}_0$), then a unique
solution to~\eqref{eqn:ivp} exists for $t \in [0, \alpha]$, where
$\alpha=\min(1, \frac{r}{M})$, and
\begin{align}
  M\equiv \max_{\vec{y} \in B_r} \|\vec{F}(\vec{y}) \|
\end{align}

In this case, $\vec{F}$ actually has continuously differentiable
component functions, which implies Lipschitz continuity, and therefore
existence and uniqueness of the solution.  Therefore, despite the lack
of continuity and monotonicity of the weak formulation of the problem,
the ODE formulation suggests there will be a unique solution.  Lack of
continuity and monotonicity also implies that one cannot prove optimal
\emph{a priori} error estimates for the weak formulation.  We will see
possible evidence of the effects of non-continuity and non-monotonicity
in the convergence results discussed in~\S\ref{sec:results}.

\section{Results\label{sec:results}}
In this section, convergence results are presented for an
implementation of the finite element formulation described
in~\S\ref{sec:fem} which is based on the libMesh
library~\cite{Peterson_2016b}.  Three representative cases are investigated:
$(\text{Re}, \alpha) = (30,15^{\circ})$, $(110,3^{\circ})$, and
$(-80,5^{\circ})$.  We also compute a ``reference'' solution using a
custom Python code~\cite{Peterson_2016a} based on the open source,
freely-available \texttt{scikits.bvp\_solver}
package~\cite{Shampine_2006,Boisvert_2012}.  This package adaptively
controls the amount of error in the numerical solution by increasing
the number of subintervals used in the calculation until a
user-defined ``tolerance'' is met.  In the present work we set the
tolerance to $10^{-14}$, which requires approximately 3200
subintervals in the most expensive case.

We remark that the \texttt{scikits.bvp\_solver} implementation is
portable, runs in under one second on a reasonably modern laptop,
and requires only about 100 lines of Python (including extensive
comments).  While nearly all authors of new solution techniques for
the Jeffery-Hamel equations compare their results to a ``reference'' solver of some type, they
typically do not provide the source code for the reference solver, and/or base it on
non-free software such as Matlab, which makes reproducing their
results difficult.  Therefore, although the code itself is
straightforward, we feel that making it readily available is, in
itself, an important contribution to the larger field.

Values of $f(\eta)$ at evenly-spaced increments in $\eta$ computed with
the $C^1$ finite element method described in \S\ref{sec:fem} are given
for comparison purposes in Table~\ref{tab:Re} for the three reference cases.
These results, which were computed using a mesh of fourth-order
Hermite elements with $h=\frac{1}{320}$, compare favorably with other
tabulated values~\cite{Motsa_2010,Inc_2013} as well as with the
reference code used in the present work.  The numerical scheme itself
performed nearly identically in each of the cases (requiring approximately
the same number of nonlinear and linear iterations) and therefore does
not appear to be particularly sensitive to the parameters $\text{Re}$
and $\alpha$.  Since the problem has such modest memory requirements,
a direct solver was actually used to precondition the linear subproblems via
PETSc's command line interface\footnote{The flag \texttt{-pc\_type lu}
was used.  The additional flag \texttt{-pc\_factor\_shift\_type nonzero}
was required to avoid a zero pivot on the finest grid with quartic elements.}.

\begin{table}[hbt]
  \centering
  \caption{Tabulated values of $f(\eta)$ for several
    different $(\text{Re}, \alpha)$ combinations
    on a mesh of 320 4th-order Hermite elements.\label{tab:Re}}
  \vspace{4pt}
  \begin{tabular}{llll}
    \toprule
    $\eta$  & $(30,15^{\circ})$       & $(110,3^{\circ})$      & $(-80,5^{\circ})$  \\
    \midrule
    0.0     &  $1.0$                & $1.0$                & $1.0$                \\
    0.1     &  $9.7312740682\e{-1}$ & $9.7923570652\e{-1}$ & $9.9596062766\e{-1}$ \\
    0.2     &  $8.9663878283\e{-1}$ & $9.1926588558\e{-1}$ & $9.8327553811\e{-1}$ \\
    0.3     &  $7.8170458993\e{-1}$ & $8.2653361228\e{-1}$ & $9.6017991246\e{-1}$ \\
    0.4     &  $6.4348113118\e{-1}$ & $7.1022118323\e{-1}$ & $9.2352159094\e{-1}$ \\
    0.5     &  $4.9758671435\e{-1}$ & $5.8049945880\e{-1}$ & $8.6845887923\e{-1}$ \\
    0.6     &  $3.5738880303\e{-1}$ & $4.4693506704\e{-1}$ & $7.8809092167\e{-1}$ \\
    0.7     &  $2.3268829344\e{-1}$ & $3.1740842757\e{-1}$ & $6.7314363566\e{-1}$ \\
    0.8     &  $1.2967274302\e{-1}$ & $1.9764109452\e{-1}$ & $5.1199108961\e{-1}$ \\
    0.9     &  $5.1642634908\e{-2}$ & $9.1230421098\e{-2}$ & $2.9155874262\e{-1}$ \\
    1.0     &  $0.0$                & $0.0$                & $0.0$                \\
    \bottomrule
  \end{tabular}
\end{table}

\begin{figure}[hbt]
  \centering
  \begin{subfigure}[t]{.45\linewidth}
    \includegraphics[width=\linewidth]{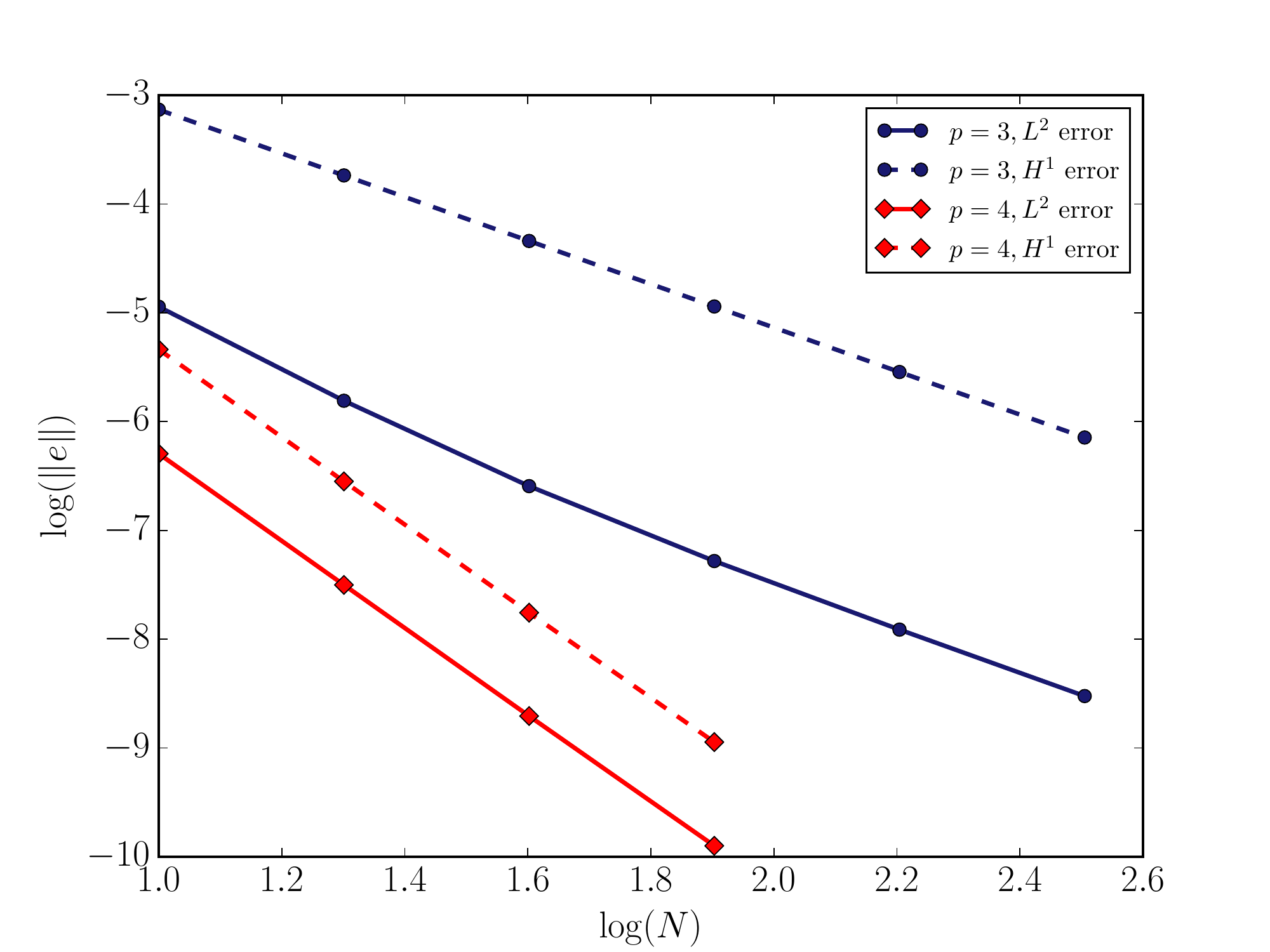}
    \caption{$(30,15^{\circ})$\label{fig:Re30_alpha15}}
  \end{subfigure}
  \begin{subfigure}[t]{.45\linewidth}
    \includegraphics[width=\linewidth]{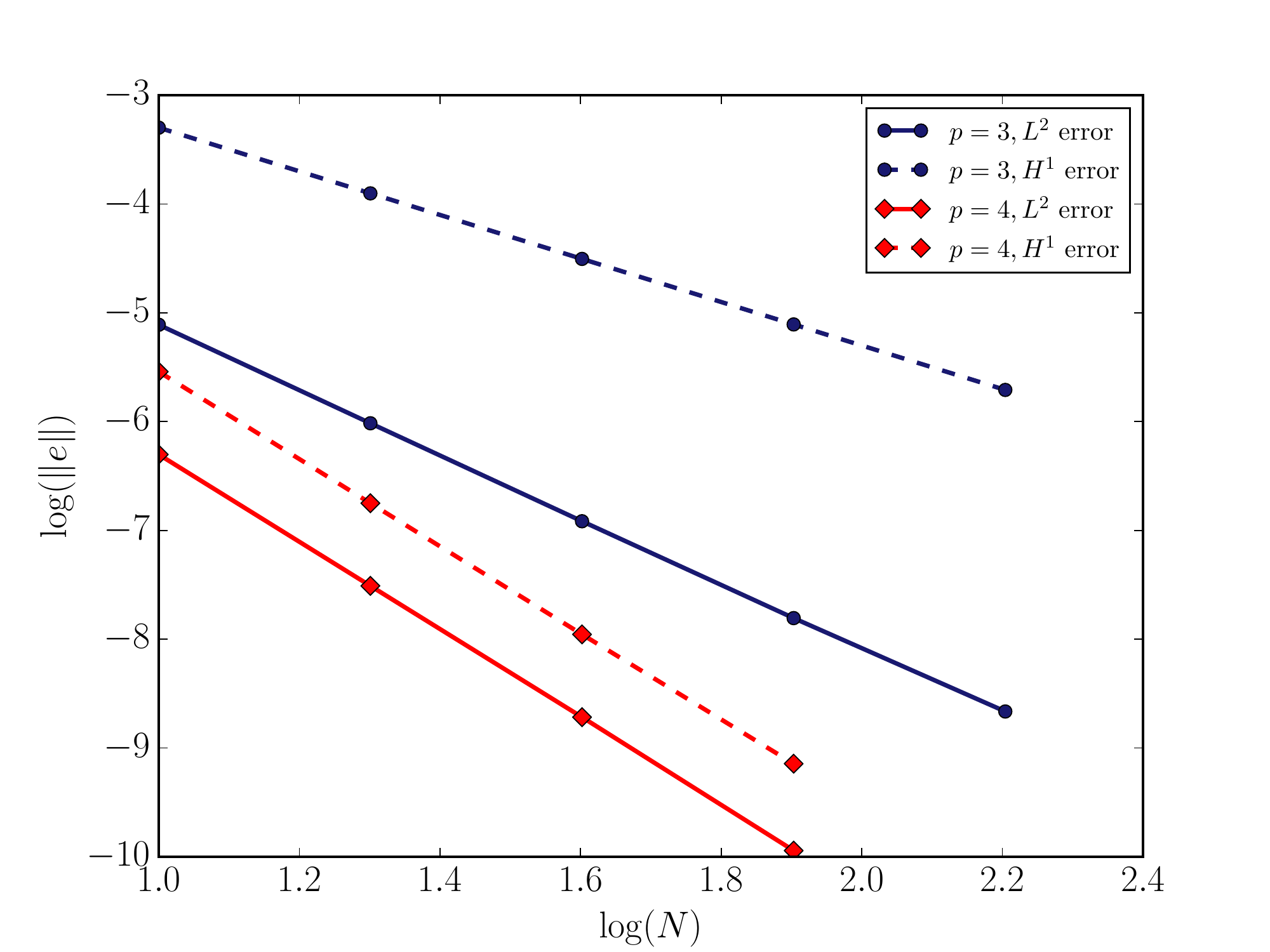}
    \caption{$(110,3^{\circ})$\label{fig:Re110_alpha3}}
  \end{subfigure}
  \\
  \begin{subfigure}[t]{.45\linewidth}
    \includegraphics[width=\linewidth]{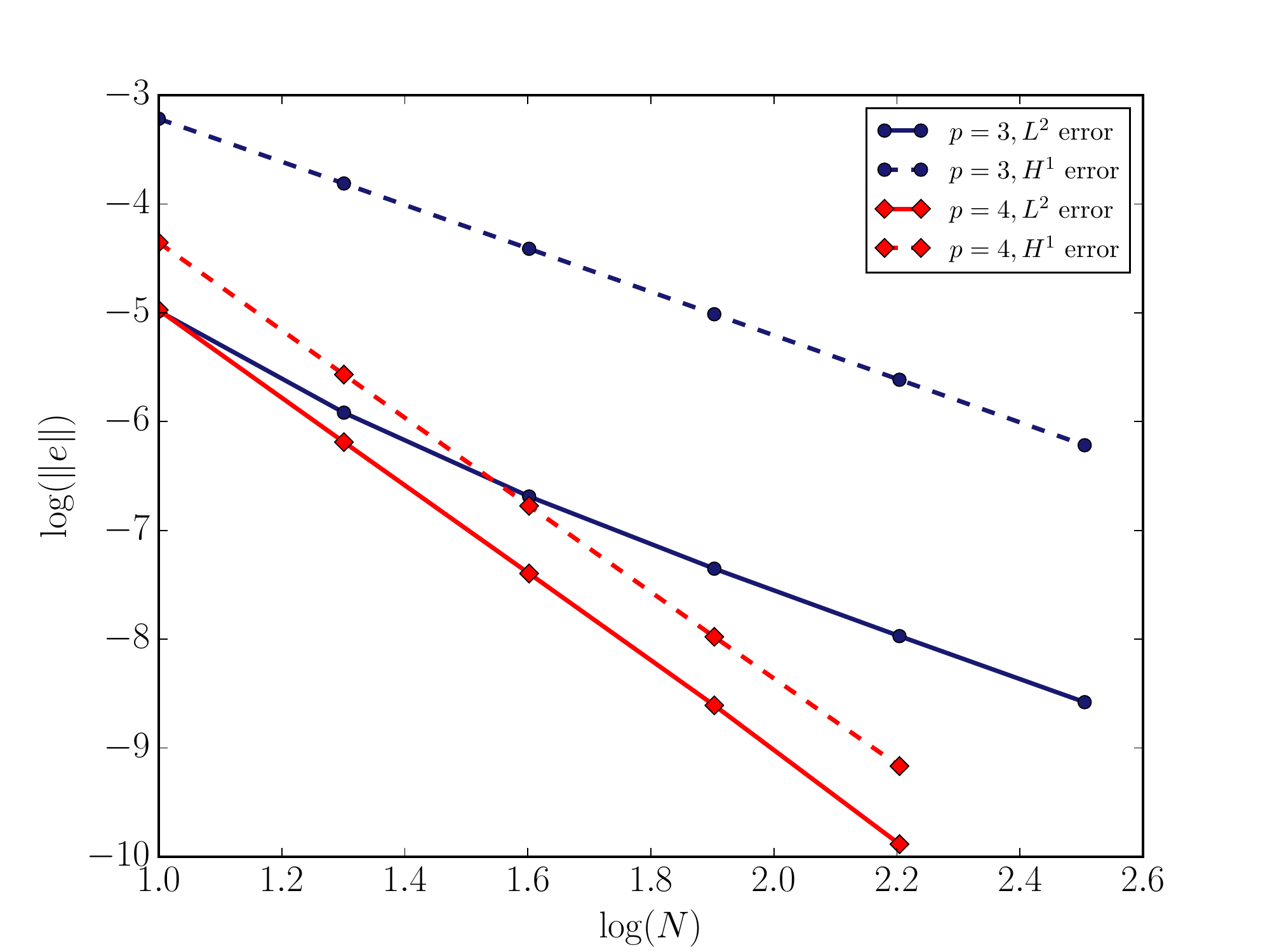}
    \caption{$(-80,5^{\circ})$\label{fig:Re-80_alpha5}}
  \end{subfigure}
  \caption{Convergence rates for representative cases $(\text{Re},
    \alpha) = (30,15^{\circ})$, $(110,3^{\circ})$, and $(-80,5^{\circ})$
    for cubic and quartic Hermite elements in the $L^2$ and $H^1$
    norms. On the $x$-axis, $N$ is equal to the number of nodes in the
    finite element mesh, and is therefore proportional to $h^{-1}$.
    \label{fig:hermite_convergence}}
\end{figure}

Convergence results for the cubic and quartic Hermite elements for the
three representative cases are given in Fig.~\ref{fig:hermite_convergence}.
Optimal convergence rates for these elements are given by \emph{a
priori} error estimation theory as $p+1$ for $\|e\|_{L^2}$ and $p$
for $\|e\|_{H^1}$, where $p$ is the polynomial degree.  Unfortunately,
we observe suboptimal convergence rates in all cases except for
the $H^1$ error on quartic elements which converges at $\mathcal{O}(h^4)$.
Furthermore, using fifth-order Hermite
elements (not shown here) also produces fourth-order accurate results
in both $L^2$ and $H^1$, reductions of 2 and 1 powers of $h$, respectively,
from the optimal rates and analogous to the observed reductions for third-order Hermite
elements.  In summary, we make the following specific observations
about the rates of convergence:
\begin{itemize}
\item The cubic elements (blue lines in
  Fig.~\ref{fig:hermite_convergence}) converge at $\mathcal{O}(h^2)$
  (see Fig.~\ref{fig:Re30_alpha15}) in some cases and $\mathcal{O}(h^3)$
  (see Fig.~\ref{fig:Re110_alpha3}) in others in the $L^2$ norm.  The
  rate of convergence therefore seems to depend on the parameters (Re,
  $\alpha$) which define the problem.
  \item The cubic elements converge at $\mathcal{O}(h^2)$ in the $H^1$ norm (blue, dashed lines) for all cases.
  \item The quartic elements converge at $\mathcal{O}(h^4)$ for both the
    $L^2$ and $H^1$ norms for all cases, with the $L^2$ norm having a smaller constant.
\end{itemize}


We currently do not have a complete explanation for the
non-optimality in the rates of convergence observed, or the
discrepancy between the rates of convergence for the even and odd
approximation orders.  For linear problems, it is well-known that
finite element formulations with non-coercive bilinear forms and
mixed formulations that don't satisfy the inf-sup condition are
unstable in the sense that they may produce reasonable results
provided that certain conditions are met ($h$ small enough, true solution
which does not excite unstable modes, etc.) but they may also
produce completely unsatisfactory (numerical oscillations, checkerboard
modes, etc.) results.  Other than the suboptimal rates of
convergence described above, we saw no evidence of unstable
modes in the third-order problem discussed here, regardless of
(Re, $\alpha$).

In non-coercive problems, stabilizing effects are sometimes achieved by
adding bubble functions to the finite element space on each element.
Examples include the so-called MINI element~\cite{Arnold_1984} for
Stokes flow in which cubic bubbles are added to linear
triangles in order to satisfy the discrete inf-sup condition, and the
addition of cubic bubble functions to stabilize convection-dominated
Galerkin discretizations of the convection-diffusion
equation~\cite{Brezzi_1992}.  It has also been observed that these
bubbles typically have no effect on either the stability or rate of
convergence for coercive problems~\cite{Franca_1994}.

In the present application, although the formulation is not unstable
in the same way that the inf-sup violating and convection-dominated
applications are unstable, adding bubbles does have a disproportionate
effect on rates of convergence in some cases ($p=3 \rightarrow 4$) but
not in others ($p=4 \rightarrow 5$), and may in some sense be said to
have a ``stabilizing'' effect on the non-coercive formulation.  More
research, especially on simplified linear model problems, is required in
order for this behavior to be fully understood.
In~\S\ref{sec:first_order_noncoercive} we briefly discuss additional numerical
results for a model problem which exhibits a similar even/odd order
convergence rate discrepancy, and might serve as the basis for further
theoretical investigations.

\subsection{First-order non-coercive problem\label{sec:first_order_noncoercive}}
\begin{figure}[hbt]
  \centering
  \includegraphics[width=.7\linewidth]{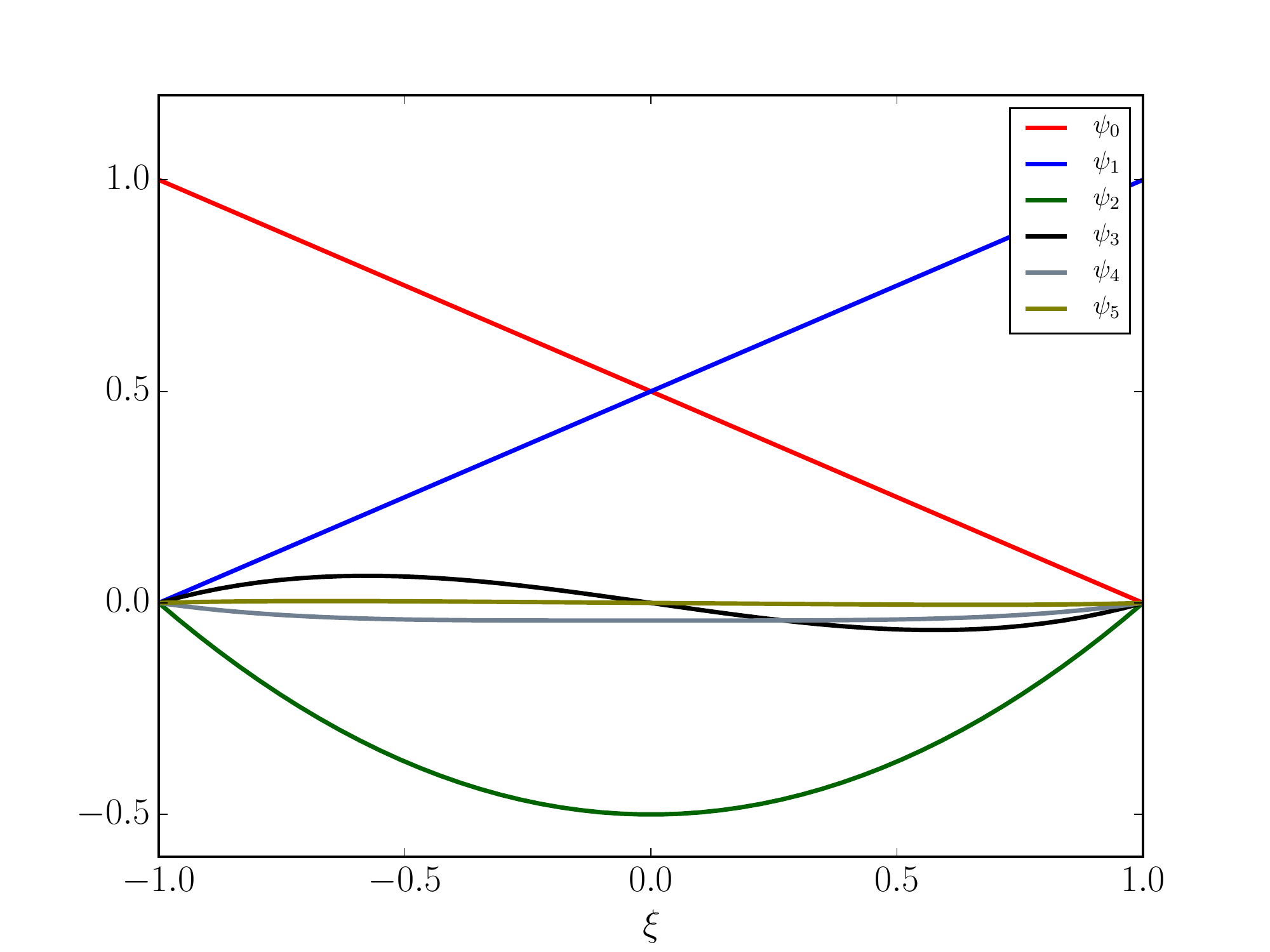}
  \caption{Hierarchic basis functions up to order $p=5$ in 1D.\label{fig:hierarchic}}
\end{figure}

\begin{figure}[hbt]
  \centering
  \begin{subfigure}[t]{.48\linewidth}
    \includegraphics[width=\linewidth]{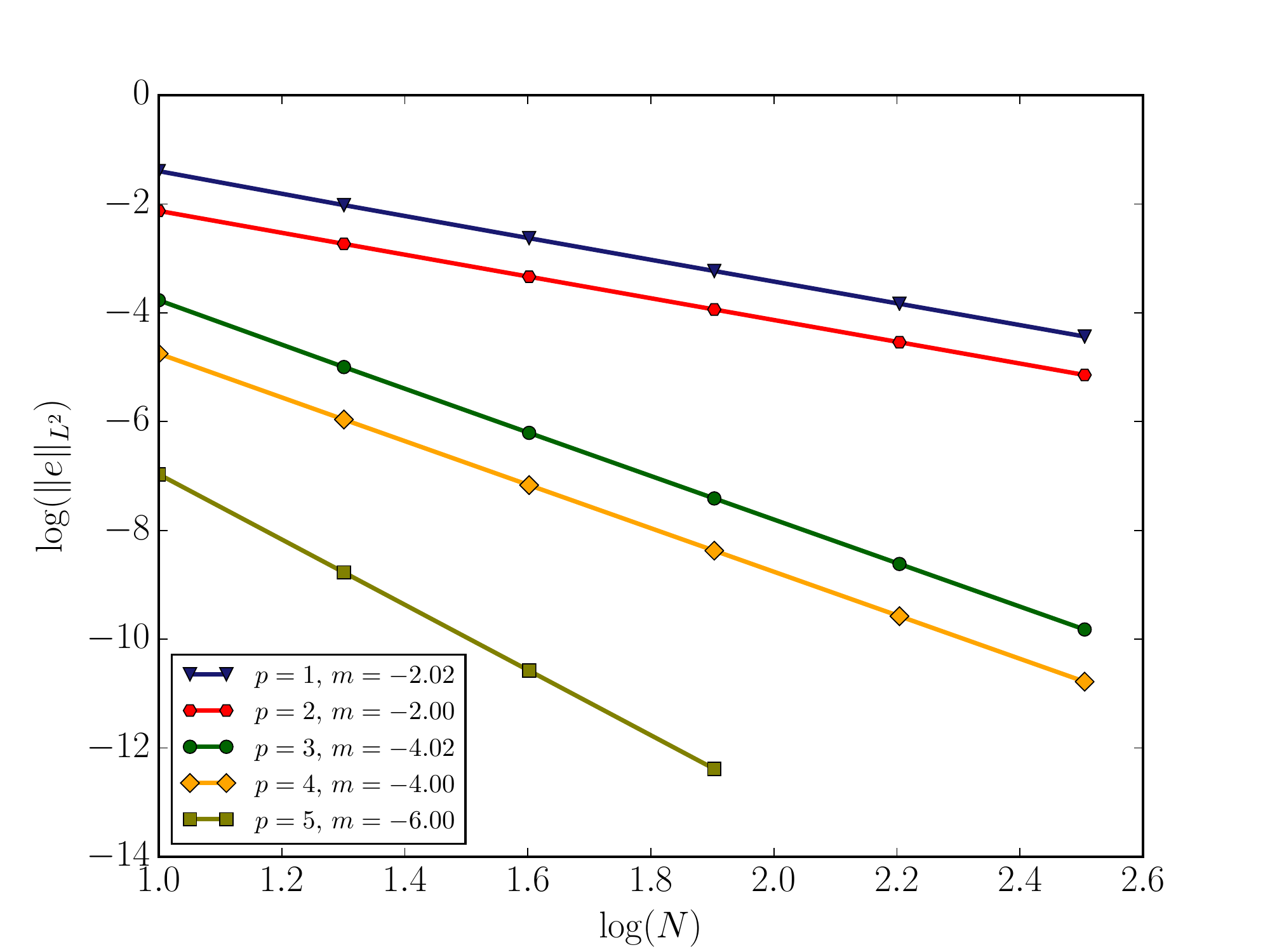}
    \caption{$L^2$ error\label{fig:hierarchic_L2}}
  \end{subfigure}
  \begin{subfigure}[t]{.48\linewidth}
    \includegraphics[width=\linewidth]{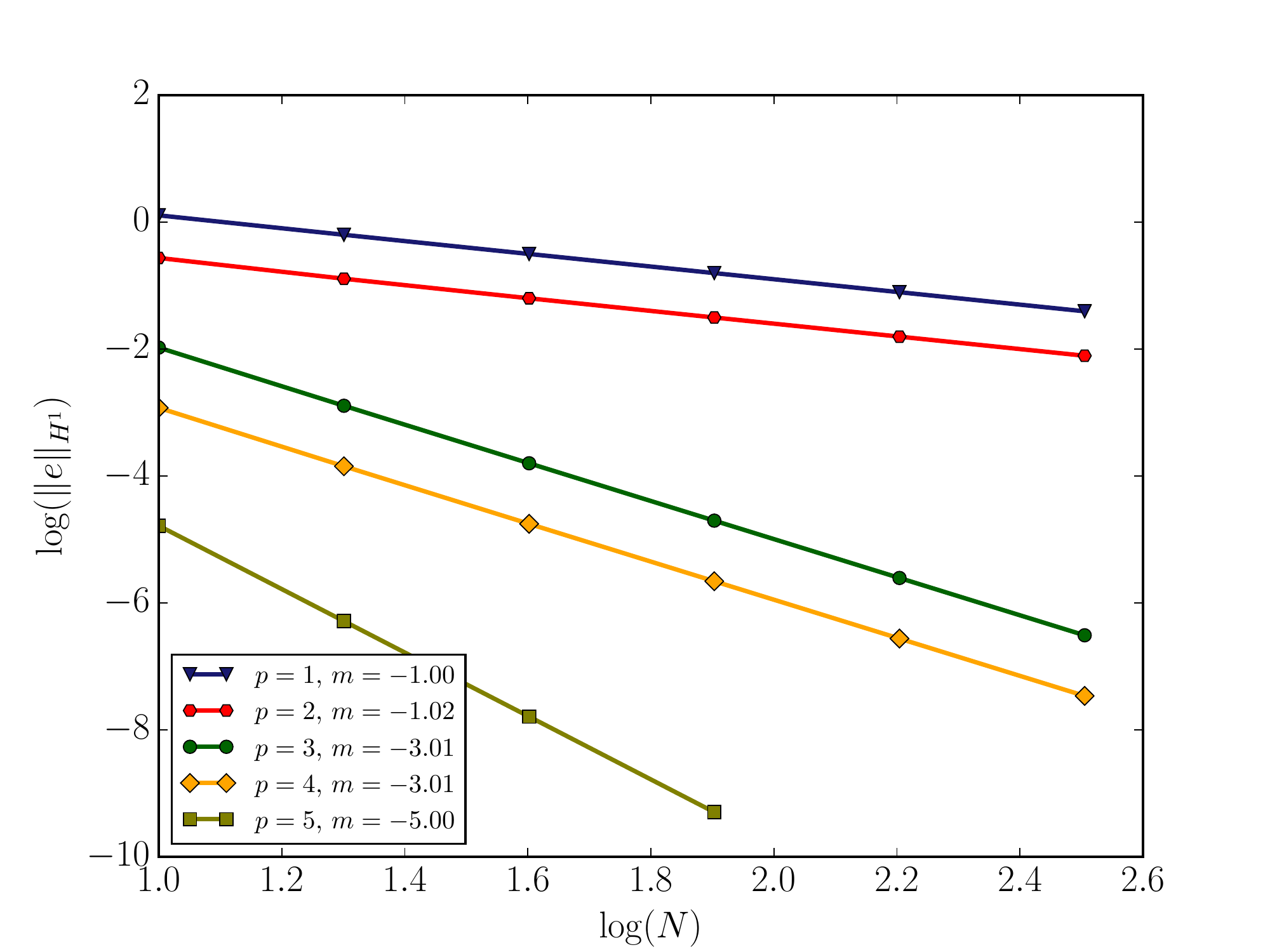}
    \caption{$H^1$ error\label{fig:hierarchic_H1}}
  \end{subfigure}
  \caption{Convergence rates, $m$, for the non-coercive first-order problem
    discretized with hierarchic finite elements of order $1 \leq p \leq 5$. The
    odd-order elements exhibit optimal convergence rates while the even-order
    elements are suboptimal by a factor of $h$ in both the $L^2$ and $H^1$
    norms of the error. Even-order elements of degree $p$ have a slightly
    better error constant (hence are more accurate) than the corresponding
    odd-order element of degree $p-1$.\label{fig:hierarchic_convergence}}
\end{figure}

To help place the different convergence rates of the even/odd order
finite element discretizations of the third-order non-coercive problem
into context, we now consider a related problem which is also
non-coercive, but is simpler to analyze due to being linear, admits a
simpler $C^0$ finite element discretization due to having only
first-order derivatives, and is trivial to compute the discretization
error for.  Specifically, we consider the ODE:
\begin{align}
  \label{eqn:ode1}
  u' + u &= g
  \\
  \label{eqn:ode1_bc}
  u(0)   &= 1
\end{align}
on $\Omega = (0,1)$, where the forcing function
\begin{align}
  g \equiv \cos \left( \frac{5 \pi x}{2} \right) -
  \frac{5 \pi}{2} \sin \left( \frac{5 \pi x}{2} \right)
\end{align}
is chosen to produce the exact solution
$u=\cos \left( \frac{5 \pi x}{2} \right)$, as may be easily verified.  The weak
formulation proceeds by multiplying~\eqref{eqn:ode1} by a test function
$v\in H^1(\Omega)$, integrating over the domain, and integrating by
parts on the leading term to obtain:
\begin{align}
  \label{eqn:ode1_weak}
  \int_0^1 \left(-uv' + uv - gv \right) \text{d}x + \left. uv \right|_0^1 = 0
\end{align}
We then incorporate the boundary conditions into the test and trial spaces by defining
\begin{align}
  V_0 &= \{v : v \in H^1(\Omega), v(0) = 0\}
  \\
  S &= \{u : u \in H^1(\Omega), u(0) = 1\}
\end{align}
and seek $u \in S$ satisfying:
\begin{align}
  \label{eqn:ode1_weak_ibp_with_bcs}
  \int_0^1 \left(-uv' + uv - gv \right) \text{d}x + u(1)v(1) = 0
  \quad \forall\, v \in V_0
\end{align}
Introducing a mesh and the finite-dimensional subspaces $V_0^h \subset V_0$, $S^h \subset S$
spanned by the basis $\{\phi_i\}, i=1..N$ leads to the discrete residual statement:
find $u_h \in S^h$ such that
$ R_i(u_h) = 0$ for $i=1..N$, where
\begin{align}
  \label{eqn:Ri_ode1}
  R_i(u_h) \equiv \int_0^1 \left(-u_h \phi_i' + u_h \phi_i - g \phi_i \right) \text{d}x
  + u_h(1)\phi_i(1)
\end{align}
The associated Jacobian contribution is independent of $u_h$ for this linear
problem, and is given by:
\begin{align}
  \label{eqn:Jij_ode1}
  J_{ij} \equiv \int_0^1 \left(-\phi_j \phi_i' + \phi_j \phi_i \right) \text{d}x + \phi_j(1)\phi_i(1)
\end{align}
We employ the ``hierarchic'' $C^0$ finite element basis functions,
which consist of the well-known linear ``hat'' functions plus bubble
functions of increasing order (see Fig.~\ref{fig:hierarchic}), to
solve~\eqref{eqn:ode1_weak_ibp_with_bcs} for different polynomial approximation orders $1
\leq p \leq 5$. The error between the finite element solution and the
known exact solution on a sequence of uniformly-refined grids is
plotted in Fig.~\ref{fig:hierarchic_convergence}.  The results show an
odd/even discrepancy in the convergence rates similar to what was seen
for the non-coercive third-order problem, however in this case it
turns out that the odd-order discretizations are optimal, while it was
the even-order discretizations which were pseudo-optimal for the
third-order problem.

Finally, we note that pursuing a standard least-squares finite element
formulation of~\eqref{eqn:ode1}--\eqref{eqn:ode1_bc} produces a
coercive bilinear form, and solutions based on hierarchic finite
elements converge to the exact solution at optimal rates in the $L^2$
and $H^1$ norms regardless of $p$.  It therefore seems reasonable that
the lack of coercivity is somehow to blame for the suboptimal
convergence rates observed in the finite element method, although we
do not attempt to develop a theoretical justification of this
observation here.

\section{Conclusions and future work\label{sec:future}}
Despite the non-optimal rates of convergence observed, the $C^1$ finite
element formulation for the third-order ODE governing Jeffery-Hamel
flow was determined to be a straightforward and numerically efficient
solution scheme.  Using $C^1$ finite elements allows the boundary
conditions to be enforced exactly within the finite element basis, and avoids the
need to iterate to determine unknown starting conditions as is
required in other boundary value problem solution methods such as the
shooting method.

The nonlinear systems of equations arising from the finite element
formulation are amenable to solution via most common linear algebra
packages, and the $C^1$ elements themselves are available in several
well-established and well-supported finite element libraries, making
the method attractive from an implementation standpoint.  The finite
element solutions have comparable accuracy to a reference boundary
value problem solution method in both the $L^2$ and $H^1$ norms of the
error on meshes of $N=160$ elements, regardless of the value of the
problem parameters (Re, $\alpha$).

The finite element formulation of the Jeffery-Hamel ODE was shown to
be non-coercive, and therefore difficult to demonstrate existence and
uniqueness for. Based on ODE existence and uniqueness arguments,
however, we do expect such a solution to exist. Additional theoretical
investigations, perhaps involving simpler, linear model problems, are
warranted to develop a more thorough explanation for the
even/odd-order discrepancy observed in the convergence rates for the
third-order problem.

A candidate problem demonstrating a similar even/odd order convergence
rate discrepancy was described and investigated numerically, but
further work is needed to develop both a convergence theory for it,
and to apply that theory to the original problem. Finally, since
third-order ODEs arise in a number of different semi-analytical and
boundary layer solutions of the Navier-Stokes equations, the $C^1$
finite element formulation developed here represents another valuable
tool in the arsenal of solution methods for such problems, and
should be easily extendable to other cases of practical interest.

\bibliographystyle{ieeetr}
\bibliography{jeffery_hamel}
\end{document}